\documentclass{ws-procs975x65}

\newcommand{\bs}[1]{\boldsymbol{#1}}
\newcommand{\pid}{\ensuremath{\mathfrak{p}}}
\newcommand{\ai}{\ensuremath{\mathfrak{o}}}
\newcommand{\aid}{\ensuremath{\mathfrak{a}}}
\newcommand{\maid}{\ensuremath{\mathfrak{m}}}
\newcommand{\ZZ}{\mathbb{Z}}
\newcommand{\QQ}{\mathbb{Q}}
\newcommand{\RR}{\mathbb{R}}
\newcommand{\CC}{\mathbb{C}}
\newcommand{\XX}{\mathbb{X}}
\newcommand{\mdim}{\operatorname{dim}^{}_{\textnormal{metr}}}
\newcommand{\hdim}{\operatorname{dim}^{}_{\textnormal{Hd}}}
\newcommand{\adim}{\operatorname{dim}^{}_{\textnormal{aff}}}
\newcommand{\ladim}{\underline{\operatorname{dim}}^{}_{\textnormal{aff}}}
\newcommand{\diam}{\operatorname{diam}}
\newcommand{\BB}{\mathcal{B}}
\newcommand{\qed}{\hfill$\Box$}

\begin{document}

\title{Iterated Function Systems in Mixed Euclidean \\ and $\pid$-adic Spaces}

\author{Bernd SING}

\address{Fakult\"{a}t f\"{u}r Mathematik\\  
Universit\"{a}t Bielefeld \\
Universit\"{a}tsstra{\ss}e 25 \\
33615 Bielefeld\\
Germany \\
E-mail: \texttt{sing@math.uni-bielefeld.de}\\
URL: \texttt{http://www.mathematik.uni-bielefeld.de/baake/sing/}
}

\maketitle

\abstracts{We investigate graph-directed iterated function systems in
  mixed Euclidean and $\pid$-adic spaces. Hausdorff
  measure and Hausdorff dimension in such spaces are defined, and an upper
  bound for the Hausdorff dimension is obtained. The relation between the Haar
  measure and the Hausdorff measure is clarified. Finally, we discus an
  example in $\RR\times\QQ^{}_2$ and calculate upper and lower bounds
  for its Hausdorff dimension.}


\section{Introduction and Setting}

The main focus of this article is the following situation: Assume that
a (finite) family $(\Omega^{}_1,\ldots,\Omega^{}_n)$ of subsets of a locally
compact Abelian group $\XX$, the topology of which is assumed to be
generated by a metric, is implicitly given as the unique solution of a
graph-directed iterated function system (GIFS). Can we
define and calculate the Hausdorff measure and Hausdorff dimension of
these sets, and determine their relation to the Haar measure in $\XX$? 

In the following, we assume that the space $\XX$ is given by
\begin{equation}\label{eq:spaceXX}
\XX=\RR^r_{}\times\CC^s_{}\times\QQ^{}_{\pid^{}_1}\times\cdots\times\QQ^{}_{\pid^{}_k}, 
\end{equation}
i.e., as a product of non-discrete locally compact fields (we shall
expand on $\QQ^{}_{\pid}$ below). We call the
number
\begin{equation}
\mdim\XX = r+2\cdot s+k
\end{equation}    
the \emph{metric dimension} of $\XX$ (also
see~Section~\ref{sec:Hausdorff_Haar}). 

\medskip

The organisation of this article is as follows: To keep everything as
self-contained as possible, we briefly review $\pid$-adic spaces in
Section~\ref{sec:padic}. In Section~\ref{sec:Hausdorff_Haar}, the
relation between Hausdorff measure and Haar measure on $\XX$ is
clarified. Iterated function systems on $\XX$ are introduced in
Section~\ref{sec:GIFS}. We define the affinity dimension for a GIFS
and show that it is an upper bound for the Hausdorff dimension of the
sets $\Omega^{}_i$. We also discus a condition for which we obtain a
lower bound for the Hausdorff dimension. In the last Section, we
explore a GIFS in $\RR\times\QQ^{}_2$.
 
\section{$\pid$-adic Spaces and their Visualisation\label{sec:padic}}

An \emph{algebraic number field} $K$ is a finite field extension of $\QQ$
lying in $\CC$, i.e., it is a simple extension of the form
$K=\QQ(\lambda)$. The \emph{integral closure} of $\ZZ$ in an algebraic
number field $K$ is called the ring of \emph{algebraic integers}
$\ai^{}_K$ of $K$. An ideal $\pid$ of the ring $\ai^{}_K$ is called
\emph{prime} if the quotient $\ai^{}_K/\pid$ is an integral domain. A
key theorem\cite{LangANT} in algebraic number theory states that every
(fractional) ideal of $\ai^{}_K$ in $K$ can be uniquely factored into
prime ideals. 

Let $K^*_{}$ be the multiplicative group of non-zero elements of
$K$. A surjective homomorphism $v: K^*_{}\to\ZZ$ with
$v(x+y)\ge\min\{v(x),v(y)\}$ (and the convention $v(0)=\infty$) is
called a \emph{valuation}\cite{Serre}. Every prime ideal $\pid$ yields a
valuation of $K$, called the $\pid$-adic valuation $v^{}_{\pid}$, and
these are all possible valuations: For $x\in K$ 
let $v^{}_{\pid}(x)=v^{}_{\pid}(x\ai^{}_K)$ (i.e., $x\ai^{}_K$ is the
(fractional) ideal generated by $x$) where a (fractional) ideal $\aid$
has the unique factorisation $\aid=\pid^{v^{}_{\pid^{}_1}(\aid)}_1
\cdots \pid^{v^{}_{\pid^{}_{\ell}}(\aid)}_{\ell}$ into prime ideals
$\pid^{}_1, \ldots, \pid^{}_{\ell}$. 

Given a $\pid$-adic valuation $v^{}_{\pid}$, one obtains an
\emph{ultrametric absolute value} (or, more precisely, a
\emph{non-Archimedean absolute value}) by
$\|x\|^{}_{\pid}=\eta^{-v^{}_{\pid}(x)}_{}$ for some $\eta>1$ (where
$\|0\|^{}_{\pid}=0$). The completion of 
$K=\QQ(\lambda)$ with respect to such a $\pid$-adic 
absolute value yields the \emph{$\pid$-adic number field}
$\QQ^{}_{\pid}$, which is a locally compact field. We note that the
completion of $\QQ$ itself w.r.t.\ the prime ideal $p\ZZ$ yields the
$p$-adic numbers $\QQ^{}_p$. 

We define the $\pid$-adic integers
$\ZZ^{}_{\pid}=\{x\in\QQ^{}_{\pid}\mathbin| \|x\|^{}_{\pid}\le 1\}$
and the related ideal $\maid^{}_{\pid}=\{x\in\QQ^{}_{\pid}\mathbin|
\|x\|^{}_{\pid}< 1\}$. Then $\ZZ^{}_{\pid}$ is a \emph{discrete
  valuation ring}, i.e., it is a principal ideal domain that has a
unique non-zero prime ideal, namely $\maid^{}_{\pid}$. Furthermore, the
\emph{residue field} $k^{}_{\pid}=\ZZ^{}_{\pid}/\maid^{}_{\pid}$ is
finite, and the choice $\eta=[\ZZ^{}_{\pid}:\maid^{}_{\pid}]$ in the
above definition of the $\pid$-adic absolute value yields the
so-called \emph{normalised $\pid$-adic absolute value} (which has nice
properties w.r.t.\ the Haar measure on $\QQ^{}_{\pid}$, see
Section~\ref{sec:Hausdorff_Haar}).

An element $\pi$ which generates $\maid^{}_{\pid}$, i.e.,
$\maid^{}_{\pid}=\pi \ZZ^{}_{\pid}$, is called a \emph{uniformizer}. By
the uniqueness of $\maid^{}_{\pid}$, the non-zero ideals of
$\ZZ^{}_{\pid}$ are given by $\pi^{m}_{}\ZZ^{}_{\pid}$
($m\in\mathbb{N}^{}_0$). If $S$ is a system of representatives of
$k^{}_{\pid}$ (including $0$ for simplicity), every element
$x\in\QQ^{}_{\pid}$ can be written uniquely as a convergent series
(w.r.t.\ the $\pid$-adic absolute value)
\begin{equation}
x=\sum_{j=m}^{\infty} s^{}_j\,\pi^j_{},
\end{equation}     
with $s^{}_j\in S$ and $m\in\ZZ$. If $x\in\ZZ^{}_{\pid}$, then one can
take $m=0$ and we simply write (with obvious meaning) $x=s^{}_0
s^{}_1 s^{}_2\ldots$. 

One can visualise $\ZZ^{}_{\pid}$ (and also $\QQ^{}_{\pid}$) as a Cantor
set\cite{Robert}. For example, if we take $\QQ^{}_2$, every
$x\in\QQ^{}_2$ can be written as $x=\sum_{j=m}^{\infty}
s^{}_j\,2^j_{}$, where $S=\{0,1\}$. Therefore, $\ZZ^{}_2$ can be
identified with the set of all $0$-$1$-sequences, i.e.,
$\ZZ^{}_2=\{0,1\}_{}^{\mathbb{N}^{}_0}$. But this is also a coding of
points in the Cantor set\footnote
{
Indeed, the Cantor set is given by $\left\{x=\frac23\cdot\sum_{j=0}^{\infty}
s^{}_j\,\left(\frac13\right)^{j}_{}\mathbin| s^{}_j\in\{0,1\}\right\}$.   
},
and points which are close in the Cantor set are also close w.r.t.\ the
$2$-adic metric (also, both the Cantor set and $\QQ^{}_2$ are totally
disconnected). In Section~\ref{sec:example}, we will use the Cantor
set to visualise sets in $\ZZ^{}_2$, where (for reasons of
representation) we take a factor of $\frac12$ instead of $\frac13$ in
the construction of the ``Cantor set'' (of course, one then obtains
the whole interval $[0,1]$).  

\section{Haar and Hausdorff Measures\label{sec:Hausdorff_Haar}}

Given an Abelian topological group $G$, a measure $\mu$ on the family
  $\BB$ of Borel sets in $G$ is called a \emph{Haar
  measure} if it satisfies the following conditions\cite{HR,Munroe,Halmos}:
\begin{itemize}
\item[\textbf{H1}] $\mu$ is a regular measure.
\item[\textbf{H2}] If $C$ is compact, $\mu(C)<\infty$.
\item[\textbf{H3}] $\mu$ is not identically zero.
\item[\textbf{H4}] $\mu$ is invariant under translations, i.e.,
  $\mu(B+t)=\mu(B)$ for all $B\in\BB$ and $t\in G$.
\end{itemize}
Haar measures are unique up to a multiplicative constant. They are
obtained by a so-called ``Method I Construction''\cite{Munroe}.  

The Haar measure on $\XX$ is the product measure of the Haar measures
of its factors. We remark that the ($1$-dimensional) Lebesgue measure
on $\RR$ and the $2$-dimensional Lebesgue measure on
$\RR_{}^2\simeq\CC$ are Haar measures. We also note that we have for
a Haar measure $\mu$ on $\RR$, resp.\ $\CC$, resp.\ $\QQ^{}_{\pid}$
\begin{itemize}
\item $\mu(\alpha B)=|\alpha|\cdot \mu(B)$ if $\alpha\in\RR$ and
  $B\subset\RR$.
\item $\mu(\alpha B)=|\alpha|^{2}_{}\cdot \mu(B)$ if $\alpha\in\CC$ and
  $B\subset\CC$.
\item $\mu(\alpha B)=\|\alpha\|^{}_{\pid}\cdot \mu(B)$ if
  $\alpha\in\QQ^{}_{\pid}$ and $B\subset\QQ^{}_{\pid}$ (where
  $\|\cdot\|^{}_{\pid}$ denotes the normalised $\pid$-adic absolute value).
\end{itemize}

On the other hand, $\XX$ is also a separable metric space, where we
take the maximum metric $d^{}_{\infty}$, i.e., for $x,y\in\XX$ with
\begin{equation}
x=(x^{}_1,\ldots,x^{}_r,x^{(1)}_{r+1}+i\cdot
x^{(2)}_{r+1},\ldots,x^{(1)}_{r+s}+i\cdot x^{(2)}_{r+s},
x^{}_{r+s+1},\ldots,x^{}_{r+s+k})
\end{equation}
(where
$x^{}_1,\ldots,x^{}_r,x^{(1)}_{r+1},\ldots,x^{(1)}_{r+s},
x^{(2)}_{r+1},\ldots,x^{(2)}_{r+s}\in\RR$, while
$x^{}_{r+s+j}\in\QQ^{}_{\pid^{}_j}$ for all $1\le j\le k$) we have 
\begin{eqnarray}
d^{}_{\infty}(x,y) &= & \max\{|x^{}_1-y^{}_1|,\ldots, |x^{}_r-y^{}_r|,
|x^{(1)}_{r+1}-y^{(1)}_{r+1}|,\ldots, |x^{(1)}_{r+s}-y^{(1)}_{r+s}|,
\nonumber \\ & & |x^{(2)}_{r+1}-y^{(2)}_{r+1}|, \ldots,
|x^{(2)}_{r+s}-y^{(2)}_{r+s}|,
\|x^{}_{r+s+1}-y^{}_{r+s+1}\|^{}_{\pid^{}_1}, \ldots,\nonumber \\ & & 
\|x^{}_{r+s+k}-y^{}_{r+s+k}\|^{}_{\pid^{}_k}\}.  
\end{eqnarray}     
Therefore, we can define the \emph{diameter} of a set $B\subset\XX$ by
$\diam(B)=\sup_{x,y\in B}d^{}_{\infty}(x,y)$ with the convention
$\diam(\varnothing)=0$. Then, the measure obtained by the so-called
``Method II Construction''\cite{Munroe,Rogers} from the set function
$\tau(B)=[\diam(B)]^{d}_{}$ is a measure and called the
\emph{$d$-dimensional Hausdorff measure} $h^{(d)}_{}$. 

Generalising~Theorem 30 in Ref.~\refcite{Rogers}, we can show that the
two measures are related as follows. 

\begin{theorem}
Let the space $\XX$ be given as in Eq.~\eqref{eq:spaceXX}, and let
$d=\mdim \XX$. Then, the $d$-dimensional Hausdorff measure $h^{(d)}_{}$
is a Haar measure. Furthermore, $h^{(d)}_{}$ equals the Haar measure
constructed as product measure where we assign measure
$1$ to the unit interval (in $\RR$) resp.\ to $\ZZ^{}_{\pid}$ (in
$\QQ^{}_{\pid}$). \qed 
\end{theorem}

As usual, it is clear that $h^{(d)}_{}(B)$ is non-increasing for a
given subset $B\subset\XX$ as $d$ increases from $0$ to
$\infty$. Furthermore, there is a unique value $\hdim B$, called the
\emph{Hausdorff dimension} of $B$, such that $h^{(d)}_{}(B)=\infty$ if
$0\le d <\hdim B$ and $h^{(d)}_{}(B)=0$ if $d>\hdim B$. 

Note that one can see from this property that Hausdorff dimension is a metric
concept rather than a topological one\cite{HW} (therefore we have chosen the
name \textit{metric} dimension; the (topological) dimension of $\XX$
is $r+2\cdot s$, because $\pid$-adic spaces $\QQ^{}_{\pid}$ are totally
disconnected).  

\section{Graph-Directed Iterated Function Systems\label{sec:GIFS}}

Let us consider the following subspace $\mathcal{L}$ of linear
mappings\footnote 
{
  One can also consider more general linear mappings\cite{Fal88}; the
  ones considered here then correspond to the case where the
  coordinate axes and the principal axes coincide.
}
from $\XX$ to $\XX$: For each $T\in\mathcal{L}$, there are numbers
$a^{}_1,\ldots, a^{}_{r+s+k}$ such that
\begin{equation}
T(x)=T((x^{}_1,\ldots,x^{}_{r+s+k}))=(a^{}_1\cdot
x^{}_1,\ldots,a^{}_{r+s+k}\cdot x^{}_{r+s+k}),
\end{equation} 
where $x^{}_1,\ldots, x^{}_r, a^{}_1,\ldots, a^{}_r\in\RR$, while
$x^{}_{r+1},\ldots, x^{}_{r+s}, a^{}_{r+1},\ldots, a^{}_{r+s}\in\CC$
and $x^{}_{r+s+j}, a^{}_{r+s+j}\in\QQ^{}_{\pid^{}_{j}}$ ($1\le j\le k$).  

We now look at the family (complex numbers $a^{}_{r+1},\ldots,
a^{}_{r+s}$ taken twice) of the $r+2\cdot s+k$ numbers
\begin{eqnarray}
& ( & |a^{}_1|,\ldots,|a^{}_r|,|a^{}_{r+1}|,|a^{}_{r+1}|,|a^{}_{r+2}|,\ldots,
|a^{}_{r+s-1}|,\nonumber \\ & &
|a^{}_{r+s}|,|a^{}_{r+s}|,
\|a^{}_{r+s+1}\|^{}_{\pid^{}_1},\ldots, \|a^{}_{r+s+k}\|^{}_{\pid^{}_k}\,),
\end{eqnarray}  
called the \emph{singular values} of $T$. We order them in descending
order $\alpha^{}_1 \ge \alpha^{}_2 \ge \ldots \ge \alpha^{}_{r+2s+k}$,
where $(\alpha^{}_1,\ldots,\alpha^{}_{r+2s+k})$ is a 
permutation of $(|a^{}_1|,\ldots,\|a^{}_{r+s+k}\|^{}_{\pid^{}_k})$. We
are only interested in maps $T\in\mathcal{L}$ which are contracting
($\alpha^{}_1<1$) and non-singular ($\alpha^{}_{r+2s+k}>0$). We
denote the subspace of non-singular and contracting maps of
$\mathcal{L}$ by $\mathcal{L}'$.

The \emph{singular value function} $\Phi^q_{}(T)$ of
$T\in\mathcal{L}'$ is defined\cite{Fal88,Fal90} for $q\ge 0$ as follows:
\begin{equation}
\Phi^q_{}(T)=\left\{ \begin{array}{l@{\quad}l}
1 & \textnormal{if } q=0 \\
\alpha^{}_1 \cdot\alpha^{}_2 \cdot\ldots\cdot \alpha^{}_{j-1}\cdot
\alpha^{q-j+1}_j & \textnormal{if } j-1<q\le j \\
(\alpha^{}_1 \cdot \alpha^{}_2\cdot \ldots \cdot
\alpha^{}_{r+2s+k})^{q/(r+2s+k)}_{} & \textnormal{if } q>r+2\cdot s+k 
\end{array} \right.
\end{equation}
Then, $\Phi^q(T)$ is continuous and strictly decreasing in
$q$. Moreover, for fixed $q$, the singular value function is
submultiplicative, i.e., $\Phi^q(T\circ U)\le\Phi^q(T)\cdot\Phi^q(U)$
for $T,U\in\mathcal{L}'$. Note that we have
$\Phi^q(T^n_{})=[\Phi^q(T)]^n_{}$. 

We now look at a \emph{graph-directed iterated function system} (GIFS)
($1\le i\le n$):
\begin{equation}
\Omega^{}_i = \bigcup_{i=1}^n \bigcup_{f^{(\ell)}_{ij}\in F^{}_{ij}}
f^{(\ell)}_{ij}(\Omega^{}_j), 
\end{equation}
where $F^{}_{ij}$ is a (finite) set of affine contracting mappings,
i.e.,
$f^{(\ell)}_{ij}(x)=T^{}_{f^{(\ell)}_{ij}}(x)+t^{}_{f^{(\ell)}_{ij}}$
with $T^{}_{f^{(\ell)}_{ij}}\in\mathcal{L}'$ and
$t^{}_{f^{(\ell)}_{ij}}\in\XX$. A GIFS 
can be visualised by a directed multi-graph
$G^{}_{(\Omega^{}_1,\ldots,\Omega^{}_n)}$, where the vertices are the sets
$\Omega^{}_i$. If $F^{}_{ij}\neq\varnothing$, we draw $|F^{}_{ij}|$
directed edges from $\Omega^{}_i$ to $\Omega^{}_j$, labelling each
edge with exactly one of the maps $f^{(\ell)}_{ij}$. We denote by
$\bs{F}$ the matrix\footnote
{
  This is also the adjacency matrix of the graph
  $G^{}_{(\Omega^{}_1,\ldots,\Omega^{}_n)}$. 
}
$\bs{F}=(|F^{}_{ij}|)^{}_{1\le i,j\le n}$ (with
the convention $|\varnothing|=0$) and by $\rho(\bs{F})$ its spectral
radius. 

We define the \emph{path space} $E^{\infty}_{}$ as the set of all infinite
paths in the graph along directed edges that start at some
vertex. Each path (and its starting point) is (uniquely, maybe after
renaming) indexed by the sequence of the edges \linebreak
$\omega=(\omega^{}_1\omega^{}_2\ldots)$ it runs along. We also
define the sets $E^{(0)}_{}=\varnothing$ (paths of length $0$), and the set
$E^{(\ell)}_{ij}$ of all paths of length $\ell$ that start at
$\Omega^{}_i$ and end at $\Omega^{}_j$ (then $\omega^{}_1\in\bigcup_{m=1}^n
F^{}_{im}$ and $\omega^{}_{\ell}\in\bigcup_{m=1}^n F^{}_{mj}$). We also set
$E^{(\ell)}_{}=\bigcup_{i=1}^n\bigcup_{j=1}^n E^{(\ell)}_{ij}$ (all
paths of length $\ell$), $E^{\textnormal{fin}}_{}=\bigcup_{\ell\ge
  0} E^{(\ell)}_{}$ (all finite paths) and
$E^{*}_{}=E^{\textnormal{fin}}_{}\cup E^{\infty}_{}$.   

For $\omega\in E^{\textnormal{fin}}_{}$ and $\varpi\in E^{*}_{}$, we
denote by $\omega\varpi$ the sequence obtained by concatenation (or
juxtaposition) if $\omega\varpi\in E^{*}_{}$. If $\omega$ is a prefix
of $\varpi$, i.e., $\varpi=\omega\ldots$, we write 
$\omega<\varpi$. By $\omega\wedge\varpi$ we denote the maximal sequence such
that both $(\omega\wedge\varpi)<\omega$ and $(\omega\wedge\varpi)<\varpi$. 
We can topologise $E^{\infty}_{}$ in a natural way using the
ultrametric $d(\omega,\varpi)=\eta^{-|\omega\wedge\varpi|}_{}$ for some
$\eta>1$. Then $E^{\infty}_{}$ is a compact space and the sets
$N(\varpi)=\{\omega\in E^{\infty}_{}\mathbin| \varpi<\omega\}$ with $\varpi\in
E^{\textnormal{fin}}_{}$ form a basis of clopen sets for $E^{\infty}_{}$. 

For $\omega=(\omega^{}_1\ldots\omega^{}_{\ell})\in
E^{\textnormal{fin}}_{}$, we define
$T^{}_{\omega}=T^{}_{\omega^{}_1}\circ\ldots\circ
T^{}_{\omega^{}_{\ell}}$ (with $T^{}_{\varnothing}(x)=x$), i.e., we
are only interested in the linear part of each map
$\omega^{}_i(x)=T^{}_{\omega^{}_i}(x)+t^{}_{\omega^{}_i}$. By the 
``Method II Construction'' with the set function
$\tau^q_{}(N(\omega))=\Phi^{q}_{}(T^{}_{\omega})$ (with
$\tau^q_{}(\varnothing)=0$), we obtain a measure $\nu^{(q)}_{}$ on
$E^{\infty}_{}$. Then, we can generalise Proposition 4.1 of
Ref.~\refcite{Fal88}.

\begin{proposition}\label{prop:affdim}
For a GIFS (with strongly connected directed graph), the following numbers
exist and are all equal:
\begin{enumerate}
\item \begin{math} \inf\{q\mathbin| \sum\limits_{\omega\in
      E^{\textnormal{fin}}_{}}
      \Phi^{q}_{}(T^{}_{\omega})<\infty\}= \sup\{q\mathbin| 
      \sum\limits_{\omega\in E^{\textnormal{fin}}_{}}
      \Phi^{q}_{}(T^{}_{\omega})=\infty\}. 
      \end{math}
\item  \begin{math} \inf\{q\mathbin| \nu^{(q)}_{}(E^{\infty}_{})=0 \}=
    \sup\{q\mathbin| \nu^{(q)}_{}(E^{\infty}_{})=\infty\}.  
      \end{math}
\item the unique $q>0$ such that\footnote
{
  As a reminder: $\rho(\bs{F})$ denotes the spectral radius of the
  matrix $\bs{F}$.
} 
\begin{equation}
\lim_{\ell\to\infty} \left(\rho\left(\left[\sum_{\omega\in
        E^{\ell}_{ij}} \Phi^{q}_{}(T^{}_{\omega})\right]^{}_{1\le
        i,j\le n}\right)\right)^{1/\ell}_{} =1.
\end{equation}
\end{enumerate}
We denote the common value by $\adim
G^{}_{(\Omega^{}_1,\ldots,\Omega^{}_n)}$ and call it the
\emph{affinity dimension} (or \emph{Falconer dimension}) of the GIFS.\qed
\end{proposition}      
   
By a covering argument, we get an upper bound for the Hausdorff
dimension of the sets $\Omega^{}_i$, compare Proposition 5.1 of
Ref.~\refcite{Fal88} and Theorem 9.12 of Ref.~\refcite{Fal90}.

\begin{proposition}
If $\nu^{(q)}_{}(E^{\infty}_{})<\infty$, then
$h^{(q)}_{}(\Omega^{}_i)<\infty$ for all $1\le i\le n$. In particular,
we have $\hdim \Omega^{}_i \le \adim
G^{}_{(\Omega^{}_1,\ldots,\Omega^{}_n)}$ for all $1\le i\le n$.\qed 
\end{proposition}

In general, it is difficult to decide whether equality holds in this
last inequality for a self-affine GIFS, although in a certain sense
equality is the generic case -- at least in $\RR^r$ (see Theorem
5.3 of Ref.~\refcite{Fal88} and Theorem 9.12 of
Ref.~\refcite{Fal90}). And contrary to the well-studied
self-similar case (where $\alpha^{}_1=\ldots = \alpha^{}_r$) in
$\RR^r$, even the \emph{open set condition}\footnote{
The OSC is satisfied if there exist disjoint non-empty
bounded open sets $(U^{}_1,\ldots,U^{}_n)$ such that  
\begin{math}
U^{}_i \supset \bigcup_{i=1}^n \bigcup_{f^{(\ell)}_{ij}\in F^{}_{ij}}
f^{(\ell)}_{ij}(U^{}_j),
\end{math}
with the unions disjoint.
}
(OSC) does not ensure the equality sign (cf.\ Ref.~\refcite{McM84} and
Examples 9.10 \& 9.11 in Ref.~\refcite{Fal90}). 

We now define a second singular value function $\Psi^q_{}(T)$ of
$T\in\mathcal{L}'$ for $q\ge 0$ as follows\cite{Fal92,Pau95}:
\begin{equation}
\Psi^q_{}(T)=\left\{ \begin{array}{l@{\quad}l}
1 & \textnormal{if } q=0 \\
\alpha^{}_{r+2s+k} \cdot\ldots\cdot \alpha^{}_{r+2s+k-j+2}\cdot\alpha^{q-j+1}_{r+2s+k-j+1} & \textnormal{if } j-1<q\le j \\
(\alpha^{}_1 \cdot \alpha^{}_2\cdot \ldots \cdot
\alpha^{}_{r+2s+k})^{q/(r+2s+k)}_{} & \textnormal{if } q>r+2\cdot s+k 
\end{array} \right.
\end{equation}
Again, $\Psi^q(T)=[\Phi^q(T^{-1})]^{-1}$ is continuous and strictly
decreasing in $q$, but supermultiplicative for fixed $q$. Just as in
Proposition~\ref{prop:affdim}, we define the \emph{lower affinity
  dimension} 
\begin{equation} 
  \ladim G^{}_{(\Omega^{}_1,\ldots,\Omega^{}_n)} =
  \inf\{q\mathbin| \sum\limits_{\omega\in E^{\textnormal{fin}}_{}}
  \Psi^{q}_{}(T^{}_{\omega})<\infty\}= \sup\{q\mathbin|
  \sum\limits_{\omega\in E^{\textnormal{fin}}_{}}
  \Psi^{q}_{}(T^{}_{\omega})=\infty\}  
\end{equation}
of the GIFS. Then, with the help of the ``mass distribution
principle'' (see Proposition 4.2 in Ref.~\refcite{Fal90}), we
obtain the following lower bound for the Hausdorff dimension of the
sets $\Omega^{}_i$, compare Proposition 2 of Ref.~\refcite{Fal92}.

\begin{proposition}\label{prop:lowaffdim}
Let $(\Omega^{}_1,\ldots,\Omega^{}_n)$ be the solution of a
(strongly connected) GIFS 
\begin{math}
\Omega^{}_i = \bigcup_{i=1}^n \bigcup_{f^{(\ell)}_{ij}\in F^{}_{ij}}
f^{(\ell)}_{ij}(\Omega^{}_j), 
\end{math}
where all unions are disjoint. If the sets
$(\Omega^{}_1,\ldots,\Omega^{}_n)$ are also pairwise disjoint, then 
$\ladim G^{}_{(\Omega^{}_1,\ldots,\Omega^{}_n)} \le \hdim \Omega^{}_i$
for all $1\le i\le n$. \qed
\end{proposition}

We remark that this disjointness condition is often easy to check in
the cases we are interested in, since $\pid$-adic spaces are totally
disconnected.  

If the linear part of all maps $f^{(\ell)}_{ij}$ is the same, i.e.,
$T=T^{}_{f^{(\ell)}_{ij}}$ for all $i,j,\ell$, we finally obtain the
following theorem. 

\begin{theorem}
For a (strongly connected) GIFS with (unique non-empty compact) solution
$(\Omega^{}_1,\ldots,\Omega^{}_n)$, where all maps $f^{(\ell)}_{ij}$ have
the same linear part $T$, the affinity dimension $\adim
G^{}_{(\Omega^{}_1,\ldots,\Omega^{}_n)}$ is given by the unique value
$q>0$ such that $\Phi^{q}_{}(T)\cdot\rho(\bs{F})=1$. The Hausdorff
dimension of the sets $\Omega^{}_i$ is bounded by the affinity
dimension of the GIFS, i.e., $\hdim \Omega^{}_i \le \adim
G^{}_{(\Omega^{}_1,\ldots,\Omega^{}_n)}$ for all $1\le i\le
n$. Furthermore, if the unions in the GIFS are disjoint and the sets  
$(\Omega^{}_1,\ldots,\Omega^{}_n)$ are pairwise disjoint, the
Hausdorff dimension of the sets $\Omega^{}_i$ is bounded from below by
the lower affinity dimension of the GIFS, i.e., $\ladim
G^{}_{(\Omega^{}_1,\ldots,\Omega^{}_n)}\le \hdim \Omega^{}_i$ for all $1\le i\le
n$, where $\ladim G^{}_{(\Omega^{}_1,\ldots,\Omega^{}_n)}$ is given by
the unique value $q>0$ such that $\Psi^{q}_{}(T)\cdot\rho(\bs{F})=1$.\qed
\end{theorem}

\section{An Example\label{sec:example}}

Our motivation for this work are so-called ``Rauzy
fractals''\cite{Rau82}, which are used to prove pure pointedness of
the dynamical system of certain $1$-dimensional sequences over a finite
alphabet, obtained by a substitution rule. ``Rauzy fractals'' yield
a geometric representation\cite{Fogg} (or so-called windows
for models sets\cite{BS04}) for such sequences.   

Here, we look at the substitution $a\mapsto aaba$, $b\mapsto aa$ 
(we obtain a two-sided infinite sequence by applying the substitution
repeatedly (we denote the zeroth position by~$|$): $a\mathbin|a\;
\mapsto\; aaba\mathbin|aaba\; \mapsto\; \ldots
aaaba\mathbin|aabaaabaaaaaba\dots$). From such a substitution, one can
obtain a GIFS (see the above literature\cite{Rau82,Fogg,BS04} and
references therein), in this case in the space $\RR\times\QQ^{}_2$:
\begin{eqnarray}
\Omega^{}_a & = & T(\Omega^{}_a) \,\cup\, T(\Omega^{}_b) \,\cup\,
T(\Omega^{}_a)+\frac12 t^{}_1 \,\cup\, T(\Omega^{}_b)+\frac12 t^{}_1 \,\cup\,
T(\Omega^{}_a)+ t^{}_2 \nonumber \\ 
\Omega^{}_b & = & T(\Omega^{}_a)+t^{}_1
\end{eqnarray}
where $T((x^{}_1,x^{}_2))=(\kappa\cdot x^{}_1,\lambda\cdot x^{}_2)$,
$t^{}_1=(\kappa,\lambda)$, $t^{}_2=(\kappa+1,\lambda+1)$,
$\kappa=\frac{3-\sqrt{17}}{2}\approx -0.562$ and
$\lambda=\frac{3+\sqrt{17}}{2}\approx 3.562$, which in the $2$-adic
expansion starts as $\lambda = 01101\ldots$. We have
$|\kappa|=\frac{2}{\lambda}$, $\|\lambda\|^{}_2=\frac12$ and
$\rho(\bs{F})=\lambda$, and therefore the affinity dimension\footnote
{
  We also have $\ladim G^{}_{(\Omega^{}_a,\Omega^{}_b)} = 2$, but the
  sets $\Omega^{}_a$ and $\Omega^{}_b$ are not disjoint. Therefore
  Proposition~\ref{prop:lowaffdim} does not apply here.
}
$\adim G^{}_{(\Omega^{}_a,\Omega^{}_b)} = 2 =\mdim
\RR\times\QQ^{}_2$. Indeed, one can show that the Haar measure of the
sets $\Omega^{}_a$ and $\Omega^{}_b$ is positive and the intersection
$\Omega^{}_a\cap\Omega^{}_b$ has Haar measure $0$.

It is more interesting to calculate the Hausdorff dimension of the
boundaries $\partial\Omega^{}_a$ and $\partial\Omega^{}_b$. For the
boundary, one can also derive a GIFS with the same contraction
$T$. This is possible, because the above GIFS for
$(\Omega^{}_a,\Omega^{}_b)$ can be dualised\cite{LW03} to obtain a point set
equation for point sets $(X^{}_a,X^{}_b)$:
\begin{eqnarray}
X^{}_a & = & T^{-1}_{}(X^{}_a) \,\cup\, 
T^{-1}_{}(X^{}_a)+ T^{-1}_{}(\frac12 t^{}_1) \,\cup\,\nonumber \\ & & 
T^{-1}_{}(X^{}_b)+ T^{-1}_{}(t^{}_1) \,\cup\, 
T^{-1}_{}(X^{}_a)+ T^{-1}_{}(t^{}_2) \\ 
X^{}_b & = & T^{-1}_{}(X_a) \,\cup\, T^{-1}_{}(X_a)+ T^{-1}_{}(\frac12
t^{}_1)\nonumber 
\end{eqnarray}
where $T^{-1}_{}((x^{}_1,x^{}_2))=(\frac1{\kappa}\cdot x^{}_1,
\frac1{\lambda}\cdot x^{}_2)$. Starting this iteration with
$X^{}_a=\{(0,0)\}=X^{}_b$, one obtains a fixed point for
$(X^{}_a,X^{}_b)$ and one can show that $J=\left(X^{}_a+\Omega^{}_a\right) \cup
\left(X^{}_b+\Omega^{}_b\right)$ is a tiling with the prototiles
$\Omega^{}_a$ and $\Omega^{}_b$ of the whole space $\RR\times\QQ^{}_2$
(for purely Euclidean spaces, this is now well
established\cite{IRpre}). With the help of this tiling $J$, one
obtains the following GIFS for the boundary: 
\begin{equation}\label{eq:IFS_boundary}
\begin{array}{lcl}
\Xi^{}_{(a,b,0)} & = & T(\Xi^{}_{(a,a,1)}) \\
\Xi^{}_{(b,a,0)} & = & T(\Xi^{}_{(a,a,-1)})+t^{}_1 \\
\Xi^{}_{(a,a,1)} & = & T(\Xi^{}_{(a,a,-1)})+t^{}_1 \,\cup\,
T(\Xi^{}_{(a,a,\frac{\lambda}2-1)}) \,\cup\,
T(\Xi^{}_{(b,a,\frac{\lambda}2-1)}) \\
\Xi^{}_{(a,a,-1)} & = & T(\Xi^{}_{(a,a,1)})  \,\cup\,
T(\Xi^{}_{(a,a,1-\frac{\lambda}2)})+\frac12 t^{}_1 \,\cup\,
T(\Xi^{}_{(a,b,1-\frac{\lambda}2)})+\frac12 t^{}_1 \\
\Xi^{}_{(a,a,\frac{\lambda}2-1)} & = & T(\Xi^{}_{(a,a,1)})+\frac12
t^{}_1 \\
\Xi^{}_{(a,a,1-\frac{\lambda}2)} & = & T(\Xi^{}_{(a,a,-1)})+t^{}_2 \\
\Xi^{}_{(a,b,1-\frac{\lambda}2)} & = &
T(\Xi^{}_{(a,a,\frac{\lambda}2-1)}) \,\cup\,
T(\Xi^{}_{(b,a,\frac{\lambda}2-1)}) \\
\Xi^{}_{(b,a,\frac{\lambda}2-1)} & = &
T(\Xi^{}_{(a,a,1-\frac{\lambda}2)})+t^{}_1 \,\cup\,
T(\Xi^{}_{(a,b,1-\frac{\lambda}2)})+t^{}_1 \\ 
\end{array}
\end{equation}     
Here, $\Xi^{}_{(a,a,1-\frac{\lambda}2)}=\Omega^{}_a \,\cap\,
\Omega^{}_a+(1-\frac{\kappa}2,1-\frac{\lambda}2)$ and similarly for the
other sets. The boundaries are therefore given by   
\begin{eqnarray}
\partial\Omega^{}_a & = & \Xi^{}_{(a,b,0)} \,\cup\, \Xi^{}_{(a,a,1)}
\,\cup\, \Xi^{}_{(a,a,-1)} \,\cup\, \Xi^{}_{(a,a,\frac{\lambda}2-1)}
\,\cup\, \Xi^{}_{(a,a,1-\frac{\lambda}2)} \,\cup\,
\Xi^{}_{(a,b,1-\frac{\lambda}2)} \nonumber\\
\partial\Omega^{}_b & = & \Xi^{}_{(b,a,0)}
\,\cup\,\Xi^{}_{(b,a,\frac{\lambda}2-1)}. 
\end{eqnarray}
To obtain a strongly connected GIFS which fulfills the disjointness
condition from the GIFS in Eq.~\eqref{eq:IFS_boundary}, we observe that
$\Xi^{}_{(a,b,0)}=\Xi^{}_{(b,a,0)}$,
$\Xi^{}_{(a,a,1)}=\Xi^{}_{(a,a,1-\frac{\lambda}2)}\cup
\Xi^{}_{(a,b1-\frac{\lambda}2)}$ and
$\Xi^{}_{(a,a,-1)}=\Xi^{}_{(a,a,\frac{\lambda}2-1)}\cup\Xi^{}_{(a,b,0)}$.
So we arrive at the GIFS
\begin{equation}\label{eq:IFS_boundary_red}
\begin{array}{lcl}
\Xi^{}_{(a,b,0)} & = & T(\Xi^{}_{(a,a,1-\frac{\lambda}2)}) \,\cup\,
T(\Xi^{}_{(a,b,1-\frac{\lambda}2)})\\ 
\Xi^{}_{(a,a,\frac{\lambda}2-1)} & = &
T(\Xi^{}_{(a,a,1-\frac{\lambda}2)})+\frac12 t^{}_1 \,\cup\,
T(\Xi^{}_{(a,b,1-\frac{\lambda}2)})+\frac12 t^{}_1\\ 
\Xi^{}_{(a,a,1-\frac{\lambda}2)} & = &
T(\Xi^{}_{(a,a,\frac{\lambda}2-1)})+t^{}_2 \,\cup\,
T(\Xi^{}_{(a,b,0)})+t^{}_2\\ 
\Xi^{}_{(a,b,1-\frac{\lambda}2)} & = &
T(\Xi^{}_{(a,a,\frac{\lambda}2-1)}) \,\cup\,
T(\Xi^{}_{(b,a,\frac{\lambda}2-1)}) \\
\Xi^{}_{(b,a,\frac{\lambda}2-1)} & = &
T(\Xi^{}_{(a,a,1-\frac{\lambda}2)})+t^{}_1 \,\cup\,
T(\Xi^{}_{(a,b,1-\frac{\lambda}2)})+t^{}_1. \\ 
\end{array}
\end{equation}     
For this GIFS, the spectral radius
$\rho(\bs{F})$ equals $2$. Consequently, we obtain $\adim
G^{}_{(\Xi^{}_{(a,b,0)},\Xi^{}_{(a,a,\frac{\lambda}2-1)},\ldots,\Xi^{}_{(b,a,\frac{\lambda}2-1)})}  
=  \frac{\log(\sqrt{17}-3)}{\log 2}+1\approx 1.1675$ and $\ladim
G^{}_{(\Xi^{}_{(a,b,0)},\Xi^{}_{(a,a,\frac{\lambda}2-1)},\ldots,\Xi^{}_{(b,a,\frac{\lambda}2-1)})}=
1$. Using the total disconnectedness of $\QQ_2$, one can show that the
disjointness condition for the sets in
Eq.~\eqref{eq:IFS_boundary_red} holds, wherefore these are 
the upper and lower bounds for the Hausdorff dimension of the
boundaries $\partial\Omega^{}_a$ and $\partial\Omega^{}_b$. We end
this article with pictures of the  
GIFS in Eq.~\eqref{eq:IFS_boundary_red} and of the sets $\Omega^{}_a$,
$\Omega^{}_b$ and their boundaries.

\begin{figure}[h]
\setlength{\unitlength}{.7cm}
\begin{picture}(15,10)
\put(10.5,8){\oval(2.5,1)}
\put(9.9,7.9){$\scriptstyle \Xi^{}_{(a,b,0)}$}

\put(4.5,5){\oval(2.5,1)}
\put(3.55,4.9){$\scriptstyle \Xi^{}_{(a,a,1-\frac{\lambda}2)}$}

\put(10.5,5){\oval(2.5,1)}
\put(9.55,4.9){$\scriptstyle \Xi^{}_{(a,a,\frac{\lambda}2-1)}$}

\put(4.5,2){\oval(2.5,1)}
\put(3.55,1.9){$\scriptstyle \Xi^{}_{(b,a,\frac{\lambda}2-1)}$}

\put(10.5,2){\oval(2.5,1)}
\put(9.55,1.9){$\scriptstyle \Xi^{}_{(a,b,1-\frac{\lambda}2)}$}

\qbezier(11.9,8)(15.9,5)(11.9,2)
\put(13.9,4.9){\vector(0,-1){0}}
\put(12.8,4.9){$\scriptstyle T(x)$}

\put(5.7,5.6){\vector(3,2){3.4}}
\put(6,7.0){$\scriptstyle T(x)+t^{}_2$}
\put(9.2,7.7){\vector(-3,-2){3.4}}
\put(7.6,6.3){$\scriptstyle T(x)$}

\put(4.5,2.6){\vector(0,1){1.8}}
\put(2.8,3.4){$\scriptstyle T(x)+t^{}_1$}

\put(10.4,4.4){\vector(0,-1){1.8}}
\put(8.5,3.4){$\scriptstyle T(x)+\frac12 t^{}_1$}
\put(10.6,2.6){\vector(0,1){1.8}}
\put(10.8,3.4){$\scriptstyle T(x)$}

\put(6,5.1){\vector(1,0){3}}
\put(7,5.3){$\scriptstyle T(x)+t^{}_2$}
\put(9,4.9){\vector(-1,0){3}}
\put(6.9,4.5){$\scriptstyle T(x)+\frac12 t^{}_1$}

\put(6,2.1){\vector(1,0){3}}
\put(7,2.3){$\scriptstyle T(x)+t^{}_1$}
\put(9,1.9){\vector(-1,0){3}}
\put(7.3,1.5){$\scriptstyle T(x)$}
\end{picture}
\caption{\hspace*{-1.05em} The directed graph $G^{}_{(\Xi^{}_{(a,b,0)},
     \Xi^{}_{(a,a,\frac{\lambda}2-1)}, 
\ldots,\Xi^{}_{(b,a,\frac{\lambda}2-1)})}$ associated to the GIFS in
Eq.~\eqref{eq:IFS_boundary_red}.}   
\end{figure}
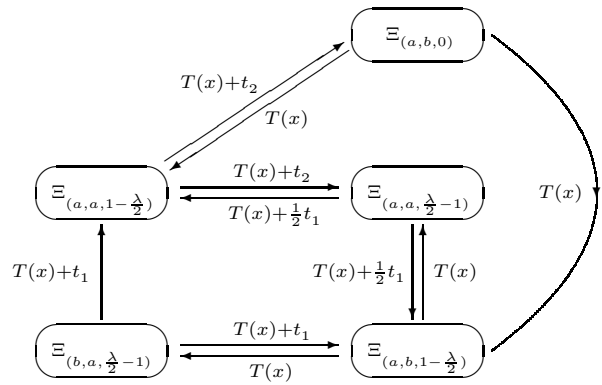

\begin{figure}[t]
\hspace*{3em}
\raisebox{14.5ex}{$\ZZ^{}_2\left\{\rule[2ex]{0em}{12ex}\right.$}
\begin{minipage}[t]{.7\textwidth} 
 \epsfxsize=\textwidth\epsfbox{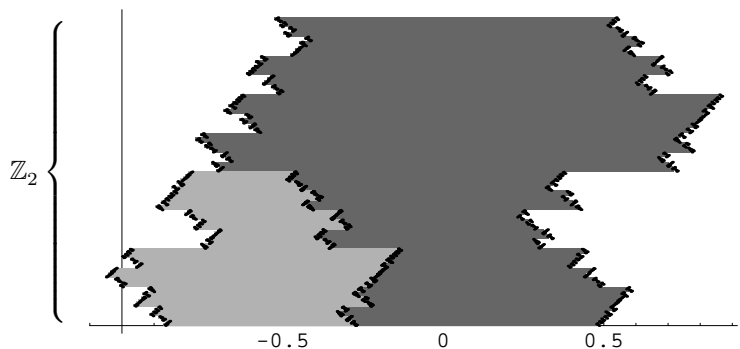} 
\end{minipage}
\caption{The sets $\Omega^{}_a$ (dark gray) and $\Omega^{}_b$ (light
  gray) and their boundaries (black) in $\RR\times\QQ^{}_2$.} 
\end{figure}

\section*{Acknowledgments}

The author thanks the referees and Michael Baake for helpful
comments. Furthermore, the author expresses his thanks to the
Cusanuswerk for financial support, as well as to the German Research
Council, Collaborative Research Centre 701.


\end{document}